\documentclass[]{article}
\usepackage{arxiv}
\usepackage[T1]{fontenc}
\usepackage{xcolor}
\usepackage[normalem]{ulem}

\usepackage{graphicx}
\usepackage{wrapfig}

\usepackage{float}

\usepackage{lipsum}
\usepackage{amsfonts}
\usepackage{graphicx}
\usepackage{epstopdf}
\usepackage{stmaryrd}
\usepackage{booktabs}
\usepackage{siunitx}
\usepackage{subcaption}

\usepackage{amsmath,amscd,latexsym}
\usepackage{amsmath}
\usepackage{mathtools}
\usepackage{bbm}
\usepackage{bm}
\usepackage[colorlinks,
    linkcolor={red!50!black},
    citecolor={blue!50!black},
    urlcolor={blue!80!black}]{hyperref}
\usepackage{nicefrac}

\usepackage{pgfplots}		
\usepgfplotslibrary{groupplots}

\newcommand{\R}{\mathbb{R}}

\newcommand{\Diffcn}{\textup{Diff}_c(\mathbb{R}^2)}
\newcommand{\Diffc}{\textup{Diff}_c(\mathbb{R}^2)}

\newcommand{\smooth}{C^{\infty}}

\newcommand{\OmegaIn}{\Omega_{\textup{in}}}
\newcommand{\OmegaOut}{\Omega_{\textup{out}}}
\newcommand{\holdall}{\Omega}

\newcommand{\id}{\textup{id}}

\usepackage[ruled,vlined,linesnumbered,algo2e]{algorithm2e} 

\SetKw{KwOr}{or}
\SetKw{KwElse}{else}
\SetKw{KwAnd}{and}
\SetKw{KwGoTo}{go to}
\title{Numerical techniques for geodesic approximation in Riemannian shape optimization}
\author{Estefania Loayza-Romero\\
University of Strathclyde, \\
Department of Mathematics and Statistics,\\  26~Richmond~St, Glasgow G1 1XH, Scotland\\
\url{estefania.loayza-romero@strath.ac.uk}
\And
Kathrin Welker \\
Helmut Schmidt University / University of the Federal Armed Forces Hamburg, \\
Faculty of Mechanical and Civil Engineering,\\ 
Holstenhofweg 85, 22043 Hamburg, Germany\\
\url{welker@hsu-hh.de}
}

\begin{document}
\maketitle              
\begin{abstract} 
Shape optimization is commonly applied in engineering to optimize shapes with respect to an objective functional relying on PDE solutions. 
In this paper, we view shape optimization as optimization on Riemannian shape manifolds. We consider so-called outer metrics on the diffeomorphism group to solve PDE-constrained shape optimization problems efficiently.
Commonly, the numerical solution of such problems relies on the Riemannian version of the steepest descent method. 
One key difference between this version and the standard method is that iterates are updated via geodesics or retractions. Due to the lack of explicit expressions for geodesics, for most of the previously proposed metrics, very limited progress has been made in this direction.
Leveraging the existence of explicit expressions for the geodesic equations associated to the outer metrics on the diffeomorphism group, we aim to study the viability of using such equations in the context of PDE-constrained shape optimization. However, solving geodesic equations is computationally challenging and often restrictive. Therefore, 
this paper discusses potential numerical approaches to simplify the numerical burden of using geodesics, making the proposed method computationally competitive with previously established methods.


\end{abstract}

\section{Introduction}

Shape optimization is commonly applied in engineering in order to optimize shapes with respect to an objective functional that relies on the solution of a partial differential equation (PDE). The PDE is required to model the underlying physical phenomenon, e.g., elastic displacements due to loadings or fluid movement due to pressure differences.  
The subject of shape optimization is covered by several fundamental monographs \cite{Allaire-2007,Delfour-Zolesio-2001,SokoZol}. 
Solving shape optimization problems is made difficult by the fact that the set of permissible shapes (called shape space) generally does not allow a vector space structure.
This fact becomes promptly clear when one tries to define the sum of two shapes and thereby realizes that this is not easily possible.
This is one of the main difficulties for the formulation of efficient methods. 
To overcome this problem, this paper sees shape optimization problems embedded in the framework of optimization on Riemannian manifolds like in \cite{geiersbach2021stochastic,Schulz2014,SchulzSiebenbornWelker,Welker2016}, i.e., the shape space is given by a Riemannian manifold. 

In the following, we are interested on solving problems of the type
\begin{equation}
\label{minproblem}
\min_{(u,y)\in \mathcal{U}\times \mathcal{Y}} \quad j(u,y)\qquad  \text{s.t. } \qquad e
(u,y)=0 ,
\end{equation}
where $j\colon \mathcal{U} \times \mathcal{Y} \rightarrow \R$ is a shape functional, $e\colon \mathcal{U} \times \mathcal{Y} \rightarrow \mathcal{W}$ is an operator, where $\mathcal{U}$ and $\mathcal{W}$ are generally Banach spaces. 
The choice of the set $\mathcal{U}$ plays a crucial role in the design of this kind of problems.
In classical shape optimization, one usually sets $\mathcal{U}$ to be a subset of the power set of some domain $\Omega$ in $\R^{n+1}$ if one considers $n$-dimensional shapes.
In this paper, we consider PDE-constrained shape optimization problems, i.e., the operator $e$ in \eqref{minproblem} represents a PDE.
As mentioned above, we view shape optimization as optimization on shape spaces, i.e., $\mathcal{U}$ is assumed to be an (infinite-dimensional) manifold endowed with a Riemannian metric $\mathcal{G}$.
A prototypical algorithm to solve PDE-constrained shape optimisation problems is the steepest descent method, for which a gradient needs to be computed. 
Since we are working on Riemannian manifolds, we use the well-known Riemannian version of this method (cf.~\cite[Alg.~4.1]{boumal2023intromanifolds}).
However, different variants of this algorithm can be devised depending on the choice of metric used to compute: first the Riemannian gradients, and second to update the iterates.
In this paper, we consider so-called outer metrics---more precisely, we consider Sobolev-type metrics on the diffeomorphism group.
Thanks to the existence the explicit expression of the geodesic equation associated to the Sobolev metrics, one could be able to use a numerical approximation of this equation in an optimization procedure to update the iterates.
In general, only approximations of the geodesics, so-called retractions, are used in the optimization methods; see, e.g., \cite{geiersbach2021stochastic,SchulzSiebenbornWelker,PryymakLoayzaWelker}. However, in \cite{herzog2024discretize} geodesics were used to update the iterates of the Riemannian steepest descent method. In this work, the advantages of using geodesics were shown, since larger step-sizes can be used without jeopardizing the quality of the meshes (since the respective metric took care of that). We will use this work, as inspiration for the following analysis. 
We will survey on efficient algorithms to numerically integrate the geodesic equation. 

The content of this paper is arranged as follows: Section \ref{sec:OptShMan} introduces a gradient-based optimization approach for shape optimization problems within the framework of Riemannian manifolds. Section \ref{sec:ShapeMan} defines the shape space and Riemannian metric considered in this paper.
We then analyze a benchmark problem from \cite{etling2020first,herzog2024discretize} in Section \ref{sec:benchmarkProblem}, covering solution existence, shape differentiability, previous numerical solution methods, and we present numerical solutions that treat shapes as elements of the diffeomorphism group and solve the problem using various outer metrics. These experiments demonstrate the advantages of outer metrics and motivate our investigation of efficient techniques for approximating the associated geodesic equations. Section \ref{sec:geodesicApproximation} concludes with potential approaches for geodesic equation approximation.

\section{Brief introduction into optimization on shape manifolds}
\label{sec:OptShMan}

We focus on solving PDE-constrained shape optimization problems, i.e., we consider the problem~\eqref{minproblem}.
A prototypical algorithm to solve such problems is the steepest descent method, for which a gradient needs to be computed. 
In this paper, we concentrate on shape manifolds as shape spaces, i.e., $\mathcal{U}$ in~\eqref{minproblem} is an (infinite-dimensional) manifold. 
There is a large number of different shape concepts, e.g., plane smooth curves \cite{Michor2007}, 
piecewise smooth curves \cite{Pryymak2023product},
surfaces in higher dimensions \cite{BauerHarmsMichor,MichorMumford2}, 
boundary contours of objects \cite{RumpfWirth2},  multiphase objects \cite{WirthRumpf}, characteristic functions of measurable sets \cite{Zolesio},  morphologies of images \cite{DroskeRumpf},  and planar triangular meshes \cite{HerzogLoayzaRomero:2020:1}. The choice of the shape space depends on the demands in a given situation. There exists no common shape space  suitable for all applications.

There are also various different types of metrics on shape spaces, e.g., inner metrics \cite{BauerHarmsMichor,Michor2007}, outer metrics \cite{BegLDDMM,Michor2007}, the Wasserstein or Monge-Kantorovic metric on the shape space of probability measures \cite{Benamou}, the Weil-Peterson metric \cite{Kushnarev}, current metrics \cite{DurrlemanPennec} and metrics based on elastic deformations \cite{FuchsJuettlerScherzerYang}. 
In general, the modeling of both, the shape space and the associated metric, is a challenging task and different approaches lead to diverse models. 

In this paper, we consider the shape manifold $\mathcal{U}$ endowed with a Riemannian metric $\mathcal{G}$. Thanks to the Riemannian metric, we not only have a classical distance measure defined by the Riemannian metric but also the gradient needed for the steepest descent method is specified. 
We want to formulate the steepest descent method to solve~\eqref{minproblem} where the shape space $\mathcal{U}$ is  a manifold endowed with a Riemannian metric $\mathcal{G}$. 
Formally, if the PDE has a (unique) solution given any choice of $u$, then the so-called control-to-state operator $S\colon \mathcal{U} \rightarrow \mathcal{Y}$, $u \mapsto y$ is well-defined. With $J(u) := j(u, S u)$ one obtains an unconstrained optimization problem $\min_{u \in \mathcal{U}} J(u)$ and this observation justifies the formulation of the algorithm for the unconstrained problem.

The shape functional $J$ is now defined on the manifold. Therefore,  for each point $u\in \mathcal{U}$, the derivative of the shape functional $J$ at $u$ is given by the push-forward $J_\ast$ associated with $J\colon \mathcal{U}\to\R$. 
With the help of the push-forward and the corresponding Riemannian metric, one can specify the Riemannian shape gradient $\nabla J(u^k)\in T_{u^k}\mathcal{U}$ with respect to $\mathcal{G}$ as $(J_\ast)_{u^k} w= \mathcal{G}_{u^k}(\nabla J(u^k),w) \, \forall w \in T_{u^k}\mathcal{U}$.
The negative Riemannian shape gradient is then used as descent direction for the objective functional $J$ in each iteration as $u^{k+1} = \exp_{u^k}\left(-s_k \nabla J(u^k)\right)$ with $s_k$ a suitable step size.
In Algorithm~\ref{alg:R-SDM}, the steepest descent method on $(\mathcal{U},\mathcal{G})$ to solve $\min_{u \in \mathcal{U}} J(u)$ is formulated.
The exponential map is defined by $\exp_{u^k}\colon T_{u^k}\mathcal{U}\to \mathcal{U},\,z\mapsto \exp_{u^k}(z)$ in the $k$-th iteration. 
Alternatively, also approximations of the exponential map, the already above-mentioned retractions, can be used to update the shape. 
Different variants of the Riemannian steepest descent method can be devised depending on the choice of Riemannian metric and retraction. 

\begin{algorithm2e}
	\caption{Steepest descent method on $(\mathcal{U},\mathcal{G})$}
	\label{alg:R-SDM}
	\DontPrintSemicolon
	\KwIn{$u^0 \in \mathcal{U} $}
    \KwData{objective functional $J\colon\mathcal{U}\to\R$ to minimize, maximum number of iterations $N_{\textup{max}}$, exponential map $\exp$ on $(\mathcal{U},\mathcal{G})$}
	Set $k \gets 0$\;
	\While{a suitable termination criterion is not satisfied or $k\le N_{\textup{max}}$}{
        Compute $\nabla J(u^k)\in T_{u^k}\mathcal{U}$ w.r.t. $\mathcal{G}$ as $(J_\ast)_u w= \mathcal{G}_u(\nabla J(u),w) \, \forall w \in T_{u^k}\mathcal{U}$ \;
		Update $u^{k+1} = \exp_{u^k}\left(-s_k \nabla J(u^k)\right)$ with $s_k$ a suitable step size\;
		Set $ k \gets k + 1 $.\;
	}
	\KwRet{$u^{k+1}$}\;
\end{algorithm2e}

\section{Riemannian shape manifold}
\label{sec:ShapeMan}

In this paper, we concentrate on one-dimensional simple, closed curves and our shape manifold is $\mathcal{U}=\Diffc/\text{Diff}_c(\R^2, \Xi)$,
where $\Xi\coloneq \{(x,y) \in \R^2 \vert \ x^2+y^2=1\} \subset \R^2$ denotes the unit circle, $ \text{Diff}_c(\R^2, \Xi):= \{ \varphi \in \Diffc \vert \varphi( \Xi)=\Xi\}$, and
$\Diffcn:= \{\varphi\in \smooth(\R^2, \R^2)|\, \varphi^{-1} \text{ exists and} \text{ is smooth}, \text{ supp}(\varphi - \id) \text{ compact}\}$ denotes the diffeomorphism group of functions equal to the identity outside a compact set. The space $\Diffcn$ is a smooth submanifold of $C^\infty(\R^2, \R^n)$ and a regular Lie group with Lie algebra given by the space of diffeomorphisms with compact support (cf. \cite[Theorem 43.1]{Kriegl1997}).
We consider the space $\Diffc$ equipped with a Sobolev-type metric of order $s$, which for two vector fields $U,V\colon \R^2\to\R^2 $, $U=U(x,y)$, $V=V(x,y)$ and $z=(x,y)$  is defined as
\begin{equation}
    \label{eq:metric1}
    H^s(U,V)= \int_{\R^2} \langle LU,V\rangle \ dz,
\end{equation}
where $L=(\operatorname{id}-A\Delta)^s$ with $A>0$ (cf. \cite{Bauer2013overview,Bauer2020sobolev}).
In the context of shape optimization, the shape space $\Diffc/\text{Diff}_c(\R^2, \Xi)$ together with these Sobolev-type metrics is firstly considered in \cite{PryymakLoayzaWelker}. 

\section{Benchmark problem}
\label{sec:benchmarkProblem}

The following example has been considered in  \cite{etling2020first,herzog2024discretize}:
\begin{equation}
\label{eq:benchmarkProb}
    \min_{u\in \mathcal{U}} \, J(u) \coloneqq \int_{\OmegaIn} y \, \textup{d} x \quad \text{s.t.} \quad \left\{\begin{array}{rl}-\Delta y & = r \text{ in } \OmegaIn \\ y &= 0 \text{ on } \partial \OmegaIn\end{array}\right.,
\end{equation}
where $\Omega = \OmegaIn \sqcup \OmegaOut$ is the hold-all domain and $u \in \mathcal{U}$ corresponds to the shape and it satisfies $u = \partial \OmegaIn$, and $r(x,y) = 2.5(x+0.4-y^2)^2 + x^2 + y^2 -1$.
We follow the reasoning in~\cite{bartels2020numerical} to understand the optimization dynamics. When minimizing the integral of $y$ over $\OmegaIn$, points $(x,y)^\top$ where $r(x,y) \leq 0$ are favorable as they reduce the value of $y$, a consequence of the maximum principle. Consequently, optimal shapes tend to resemble the $0$-th sublevel set of $r$.

The existence of solutions and shape differentiability of this problem has been studied in~\cite{etling2020first}, and thus without proof we present the weak formulation of its shape derivative in the direction of a sufficient continuous vector field~$\bm{W}$:
\begin{equation}
\label{eq:shapederivative}
    \textup{d} J(u)[\bm{W}]=  \int_{\OmegaIn} y \operatorname{div}(\bm{W}) \, \textup{d} x
    - \int_{\OmegaIn} \operatorname{div}(r\bm{W})p \, \textup{d} x 
    + \int_{\OmegaIn} (\nabla y)^\top [\operatorname{div}(\bm{W}) 
- D\bm{W} - D\bm{W}^\top]\nabla p \, \textup{d} x   ,
\end{equation}
where the adjoint state $p$ is the unique solution of the problem
\begin{equation}
\left\{
\begin{array}{rll}
    -\Delta p &= -1  &\quad\text{in} \, \OmegaIn, \\ 
     p &= 0 &\quad \text{on} \, \partial \OmegaIn.
     \end{array}
    \right.
\end{equation}

The discretized version of problem~\eqref{eq:benchmarkProb} is obtained by representing the hold-all domain and the corresponding shape $u$ in terms of a triangular mesh. We will refer to the domain covered by the mesh $\Omega_h$. We will use the finite element method to discretize the state and adjoint equations in~\eqref{eq:benchmarkProb}. Thus, $\mathcal{S}^1(\Omega_h)$ denotes the finite element space of piecewise linear, globally continuous functions. Moreover, let us denote by $\mathcal{S}^1_0(\Omega_h)$ the subspace of functions with zero Dirichlet boundary conditions. 

As a motivation, we start mentioning that this problem has been previously addressed using various methods. Our primary comparison reference is~\cite{herzog2024discretize}, who employed a discretize-then-optimize approach for problem \eqref{eq:benchmarkProb}. They proposed a complete Riemannian metric for planar triangular meshes (a finite-dimensional manifold introduced in \cite{HerzogLoayzaRomero:2020:1}), and derived its associated geodesic equation. This enabled them to implement a variant of the Riemannian steepest descent method with geodesic updates, termed the Complete-Complete (Comp-Comp) variant. They also explored a second Riemannian metric---the Elastic-Euclidean (Elas-Euc)---which represents the discrete counterpart of the Steklov-Poincaré (SP)-metric from \cite{SchulzSiebenbornWelker}, where ``Euclidean'' refers to updating iterates along Euclidean geodesics (retraction). 
    
Their results showed that while the variant Elas-Euc required $150$ iterations to reach a value of the objective function of $-0.0912$ in a total of $\SI{1.726}{\second}$. 
Notably, the variant Complete-Complete, i.e., the one using the complete metric to compute the shape gradients and updating the iterates using the associated geodesic equation only required $15$ iterations to reach a value of the objective function of $-0.0913$.
The authors attribute this behavior to the fact that thanks to the completeness of the metric the algorithm can take longer step sizes without jeopardizing the quality of the mesh. 
The main drawback of this method is the total computation time is about $\SI{89.761}{\second}$. 

\begin{figure}[t]
\captionsetup[subfigure]{labelformat=empty}
    \begin{center}
        \begin{subfigure}{0.24\textwidth}
            \includegraphics[width=\textwidth]{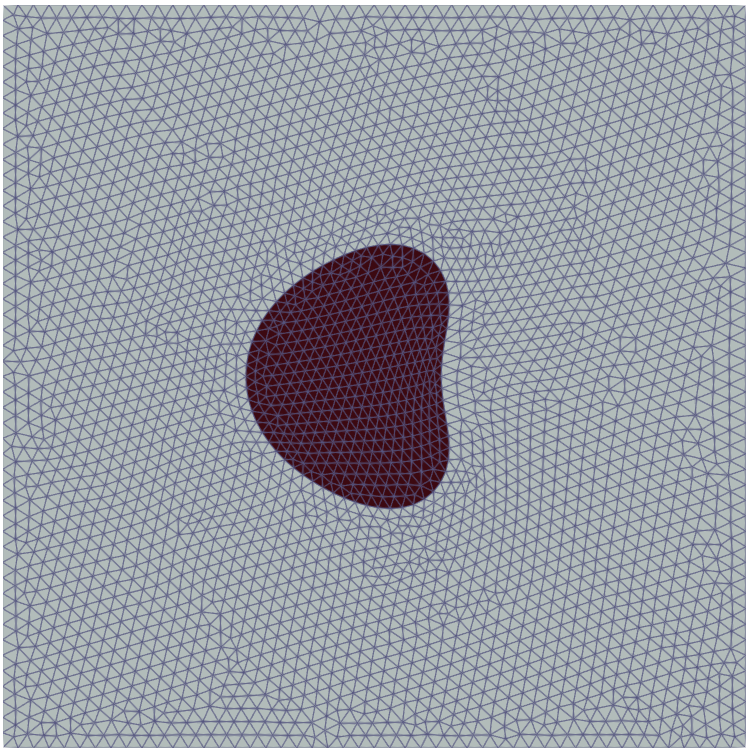}
            \caption{SP-metric: $\OmegaIn^{158}$}
        \end{subfigure}
        \begin{subfigure}{0.24\textwidth}
            \includegraphics[width=\textwidth]{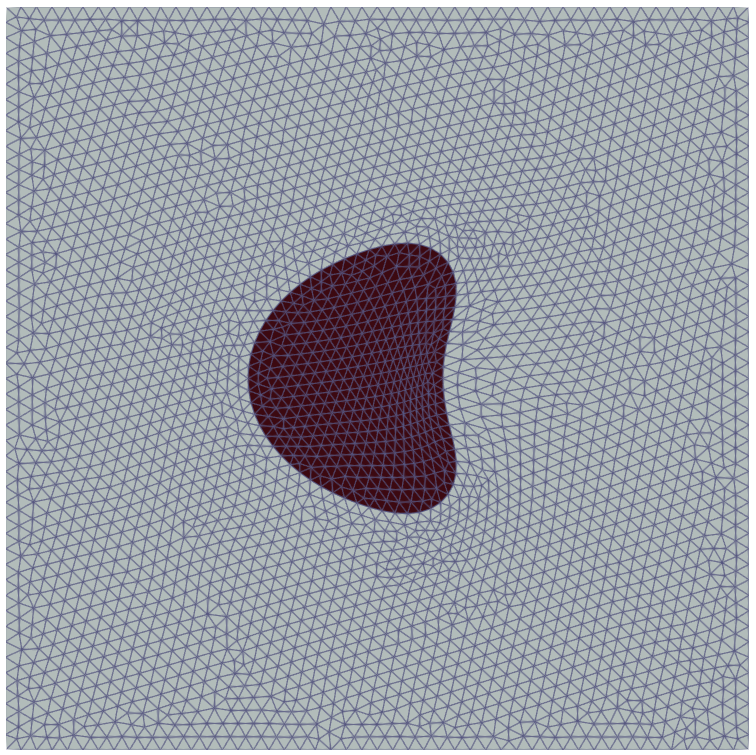}
             \caption{$H^2$-metric: $\OmegaIn^{40}$}
        \end{subfigure}
        \begin{subfigure}{0.24\textwidth}
            \includegraphics[width=\textwidth]{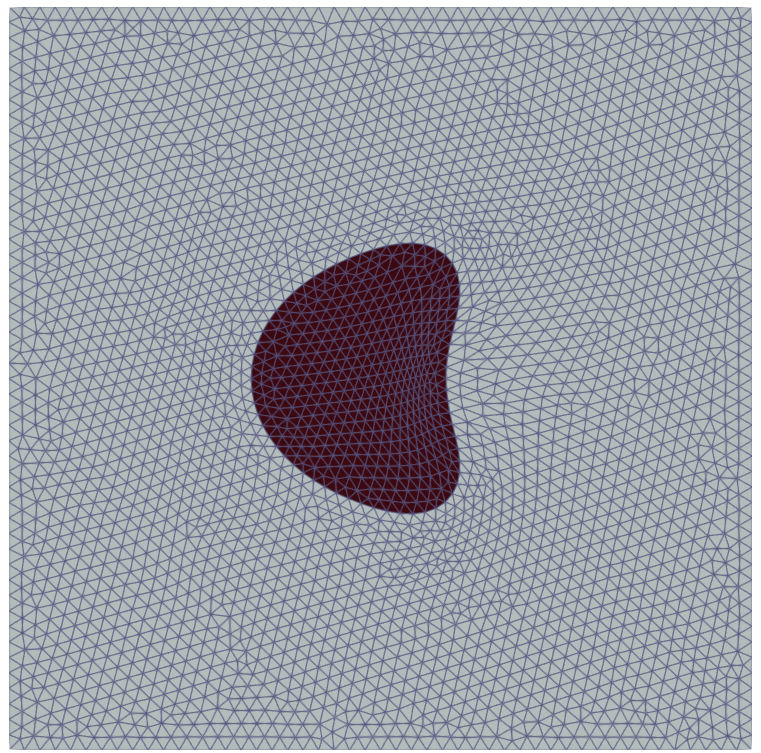}
             \caption{$H^3$-metric: $\OmegaIn^{36}$}
        \end{subfigure}
        \begin{subfigure}{0.24\textwidth}
            \includegraphics[width=\textwidth]{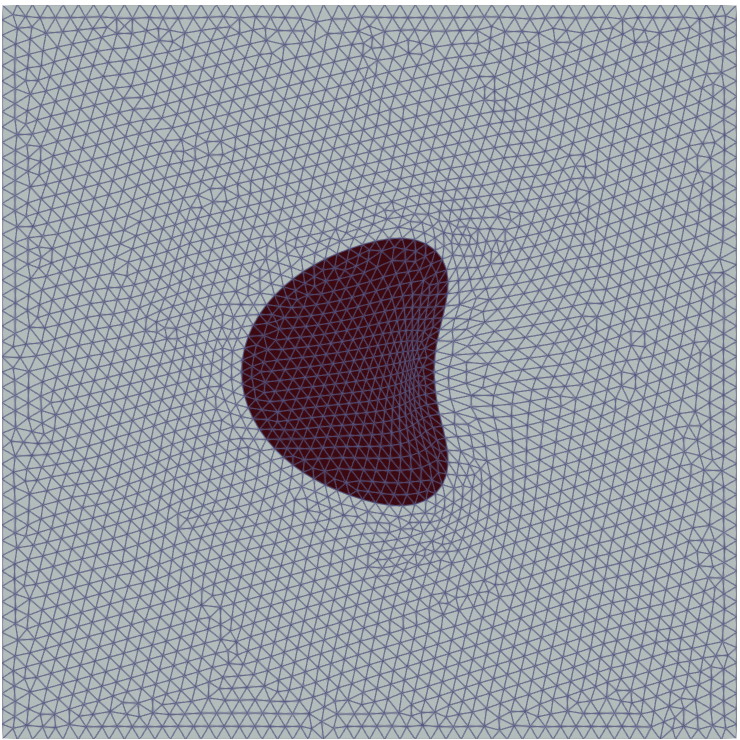}
             \caption{$H^4$-metric: $\OmegaIn^{31}$}
        \end{subfigure}
        \begin{subfigure}{0.24\textwidth}
            \includegraphics[width=\textwidth]{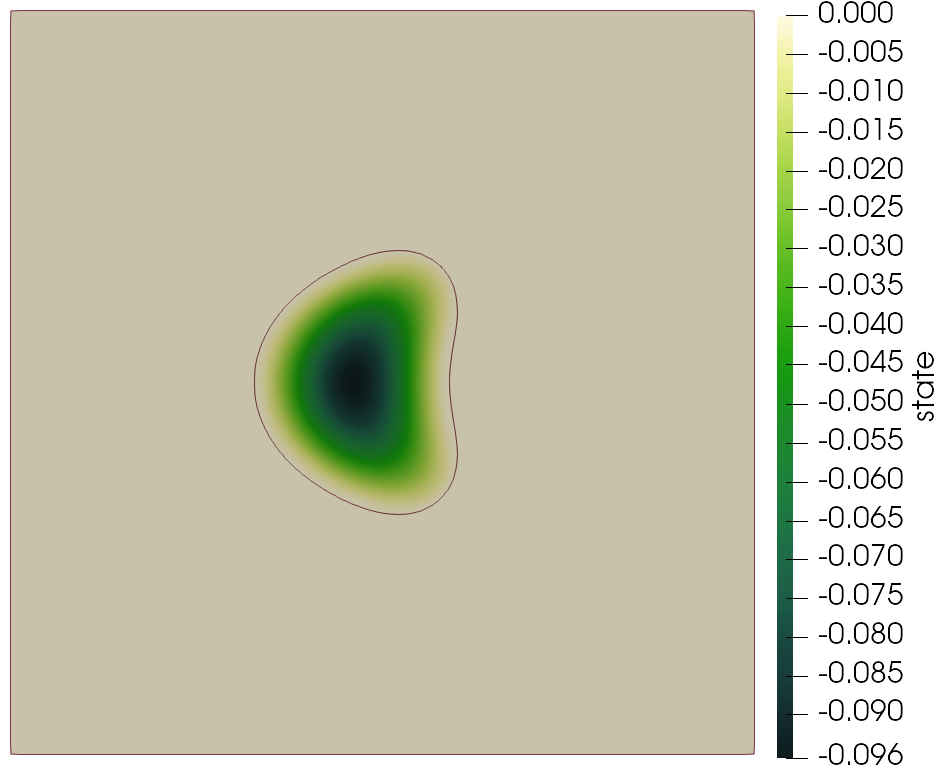}
             \caption{SP-metric: $y\left(\OmegaIn^{158}\right)$}
        \end{subfigure}
        \begin{subfigure}{0.24\textwidth}
            \includegraphics[width=\textwidth]{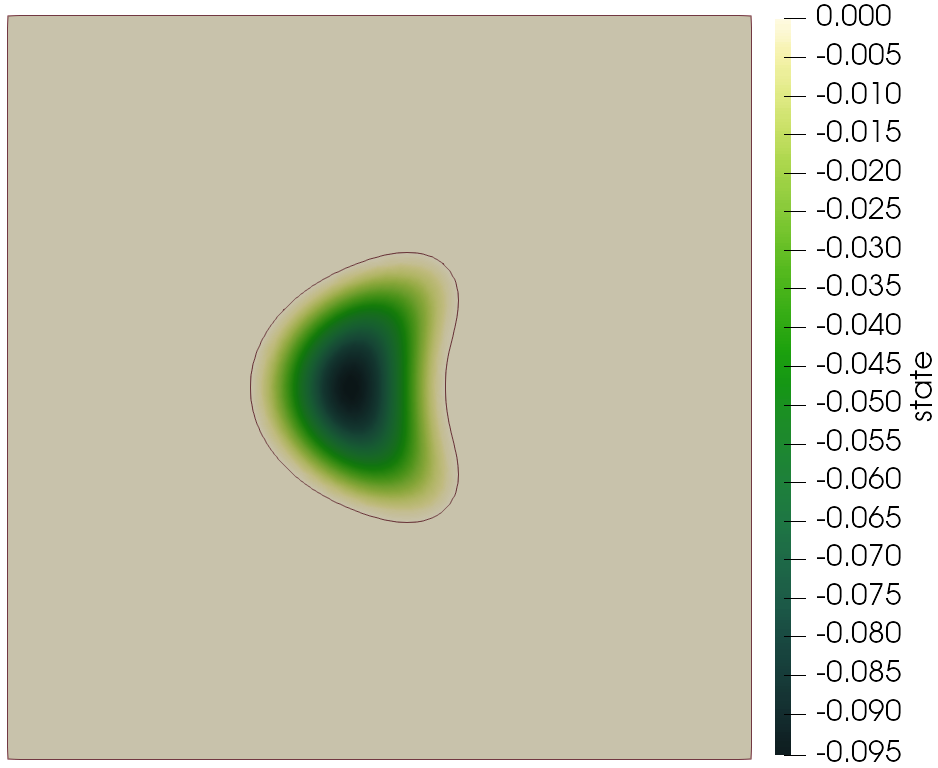}
            \caption{$H^2$-metric: $y\left(\OmegaIn^{40}\right)$}
        \end{subfigure}
        \begin{subfigure}{0.24\textwidth}
            \includegraphics[width=\textwidth]{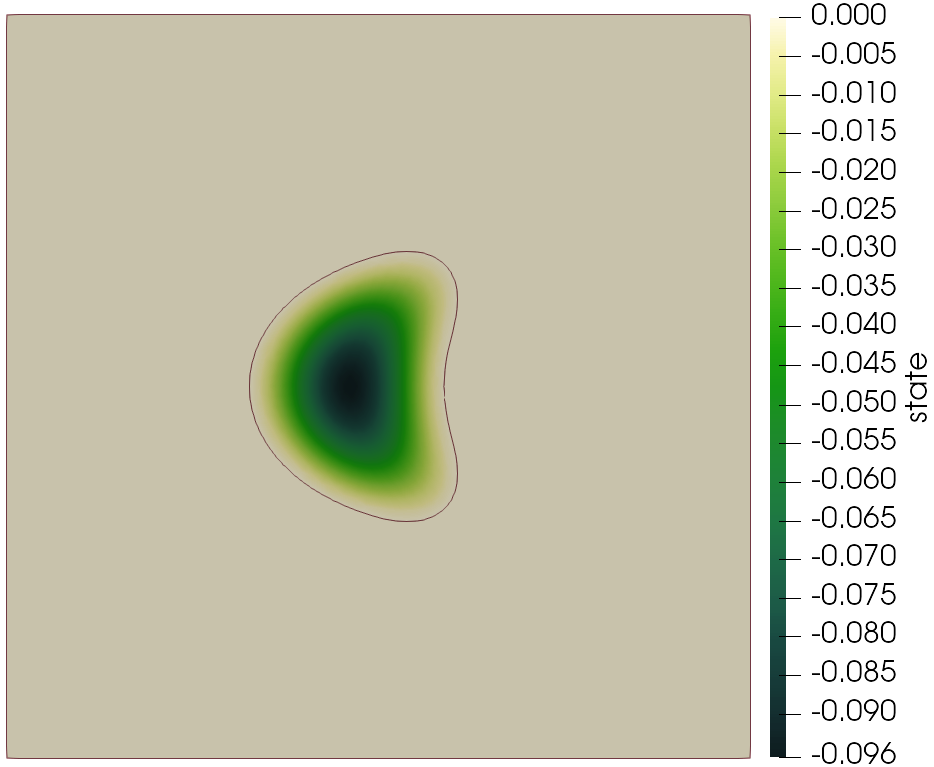}
            \caption{$H^3$-metric: $y\left(\OmegaIn^{36}\right)$}
        \end{subfigure}
        \begin{subfigure}{0.24\textwidth}
            \includegraphics[width=\textwidth]{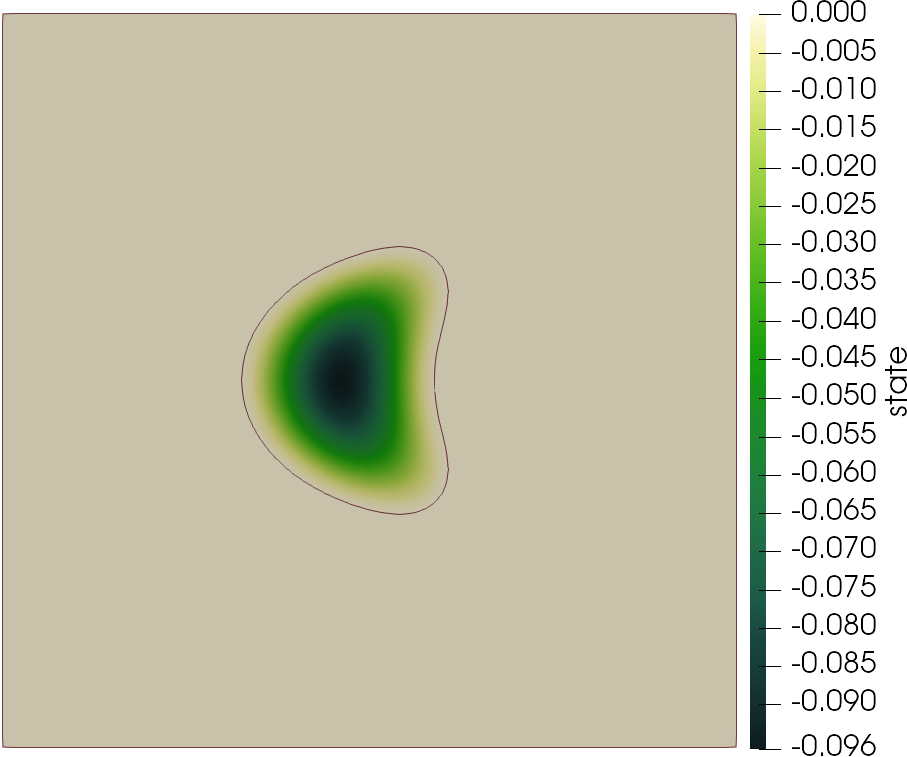}
            \caption{$H^4$-metric: $y\left(\OmegaIn^{31}\right)$}
        \end{subfigure}
    \end{center}
    \caption{Final shapes (upper row) and associated state variables (lower row) for the experiment described in Section~\ref{sec:benchmarkProblem}.}
    \label{fig:shapeNewRes}
\end{figure}

On the other hand, recent research has highlighted outer metrics over the space of diffeomorphisms (introduced in \cite{PryymakLoayzaWelker} and described in equation~\eqref{eq:metric1}) as a promising alternative to elasticity metrics. These outer metrics offer several advantages: they can be easily implemented in FEniCS, making them accessible for various applications, and they appear to maintain mesh quality more effectively than the SP-metric, which is an inner metric. To validate these claims, we numerically solve problem~\eqref{eq:benchmarkProb} using the outer metrics given in the previous section.
A key distinction from previously proposed methods is that outer metrics require defining a hold-all domain where the metric operates. For our implementation, we set this domain as the box $[-3,3]\times[-3,3]$ and use an initial shape of a circle centered at $(0,0)^\top$ with radius $1$. The initial mesh consists of $3435$ nodes and  $6932$ elements. We solve the problem using $H^2$-, $H^3$-, and $H^4$-metrics. Figure~\ref{fig:shapeNewRes} displays the different solutions obtained through various versions of the Riemannian steepest descent method.

As stopping criterion we use $\max_{m = 1,\ldots, 10} J(u^{k-m}) - J(u^k) < 10^{-4}$ and for all experiments we set the step-size to be fixed and equal to $1$. The SP-metric uses the parameters $\lambda = 0 $ and $\mu$ obtained by solving a Poisson's equation with $\mu_{\textup{min}} = 1$ at the outer boundary and $\mu_{\textup{max}} = 5$ at the moving boundary. 
The parameters used for the $H^s$-metrics have been described in Table~\ref{tab:NumExp}. 

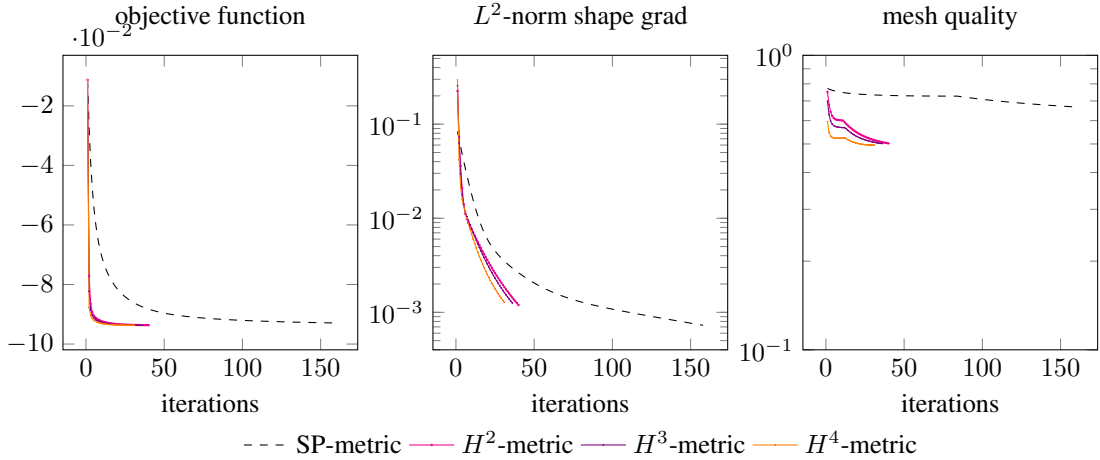
\begin{figure}[t]
\pgfplotsset{/pgfplots/group/.cd, horizontal sep=1cm}
\begin{center}
\begin{tikzpicture}
\pgfplotsset{compat=1.14, width=\textwidth}
\begin{groupplot}[group style={
        {group size=3 by 1}, xlabels at=edge bottom}, height=5.5cm, width=5.5cm,
        legend style={transpose legend, legend columns=4, draw=none }, 
        xlabel=iterations]
\nextgroupplot[title= objective function, legend to name=grouplegend, y tick label style={scaled ticks=base 10:2},scaled x ticks=false]
\addplot[color=black, dashed] table [x=iter,y=objective, col sep=comma]{Data/history_s_0.txt};\label{pgfplots:SP}
\addlegendentry{SP-metric}
\addplot[color=magenta, mark=*, mark size = 0.25pt] table [x=iter,y=objective, col sep=comma]{Data/history_s_2.txt};\label{pgfplots:H2}
\addlegendentry{$H^2$-metric}
\addplot[color=violet, mark=x, mark size = 0.25pt] table [x=iter,y=objective, col sep=comma]{Data/history_s_3.txt};\label{pgfplots:H3}
\addlegendentry{$H^3$-metric}
\addplot[color=orange, mark =star, mark size = 0.25pt] table [x=iter,y=objective, col sep=comma]{Data/history_s_4.txt};\label{pgfplots:H4}
\addlegendentry{$H^4$-metric}

\nextgroupplot[title =$L^2$-norm shape grad, ymode=log]
\addplot[color=black, dashed] table [x=iter,y=norm_felas, col sep=comma]{Data/history_s_0.txt};\label{pgfplots:SP}
\addplot[color=magenta, mark=*, mark size = 0.25pt] table [x=iter,y=norm_felas, col sep=comma]{Data/history_s_2.txt};\label{pgfplots:H2}
\addplot[color=violet, mark=x, mark size = 0.25pt] table [x=iter,y=norm_felas, col sep=comma]{Data/history_s_3.txt};\label{pgfplots:H3}
\addplot[color=orange, mark =star, mark size = 0.25pt] table [x=iter,y=norm_felas, col sep=comma]{Data/history_s_4.txt};\label{pgfplots:H4}  

\nextgroupplot[title = mesh quality, ymode=log, ymin= 1e-1, ymax = 1]
\addplot[color=black, dashed] table [x=iter,y=msh_quality, col sep=comma]{Data/history_s_0.txt};
\addplot[color=magenta, mark=*, mark size = 0.25pt] table [x=iter,y=msh_quality, col sep=comma]{Data/history_s_2.txt};
\addplot[color=violet, mark=x, mark size = 0.25pt] table [x=iter,y=msh_quality, col sep=comma]{Data/history_s_3.txt};
\addplot[color=orange, mark =star, mark size = 0.25pt] table [x=iter,y=msh_quality, col sep=comma]{Data/history_s_4.txt};
\end{groupplot}

\node[yshift=-35pt] at ($(group c2r1.south)!0.5!(group c2r1.south)$) {\ref{pgfplots:SP} SP-metric  \ref{pgfplots:H2} $H^2$-metric  \ref{pgfplots:H3} $H^3$-metric  \ref{pgfplots:H4} $H^4$-metric};
\end{tikzpicture}
\caption{Objective function, $L^2$-norm of shape gradients and mesh quality behavior along the iterations for the different variants of the method for the experiment described in Section~\ref{sec:benchmarkProblem}.}
\label{fig:valuesNewRes}
\end{center}
\end{figure}

\begin{table}[h]
\begin{center}
\begin{tabular}{c|c|crrr}
\toprule
metric  & $A$ & $\bar{k}$ &\multicolumn{1}{c}{ $J(u^{\bar{k}})$} &\multicolumn{1}{c}{ $\|V^{\bar{k}}\|_{L^2}$}& \multicolumn{1}{c}{$\varphi(\holdall_h^{\bar{k}})$} \\ 
\midrule
SP & \multicolumn{1}{c|}{--} & 158 & $-0.092910$ &  $7.28\times 10^{-3}$ & $6.68\times 10^{-1}$ \\
\midrule
$H^2$  & 0.09  & 40 & $-0.093664$ & $1.19\times 10^{-3} $& $5.02\times 10^{-1}$\\
$H^3$  & 0.04  & 36 & $-0.093684$ & $1.25\times 10^{-2}$ & $5.01\times 10^{-1}$\\
$H^4$  & 0.02 & 31 & $-0.093699$ & $1.27\times 10^{-2}$ & $4.95\times 10^{-1}$
\end{tabular}
\end{center}
\caption{Summary of the experiment described in Section~\ref{sec:benchmarkProblem}.}
\label{tab:NumExp}
\end{table}
Moreover, we show in Figure~\ref{fig:valuesNewRes}, the values of the objective function and mesh quality obtained for each metric. The experiments show the use of outer metrics reduce the amount of iterations required to achieve the same lower values of the objective function compared with the SP-metric. Moreover, since now we are allowed to use finer meshes, one can observe that the value of the objective function is lower that the one obtained in~\cite{herzog2024discretize}. 

Finally, we draw attention to the fact that the quality of the meshes remains within acceptable values even for the SP-metric. This is computed as the minimum of the radius ratio of cells from the triangulation. The radius radio is given by the quotient of the inradius and circumradius times a normalization factor depending on the geometric dimension of the problem, given by 
\begin{equation}
    \label{eq:mesh_quality}
\varphi(\holdall_h^k) = \min_{T \in \Sigma_h^k} 2\frac{\textup{inrad}(T)}{\textup{circumrad}(T)} 
\end{equation}
where $\Sigma_h^k$ denotes the triangulation associated to the mesh describing $\holdall_h^k$, and $T$ denotes a triangle from the triangulation. This quality measure has range zero to one where close-to zero values indicate the presence of at least one degenerate element.

The main aim of this experiment is to show that even without the use of geodesics, outer metrics have overall better performance than the SP-metric. Which in turn also justifies why we are interested on using their associated geodesic equation to update the iterates for this variant of the Riemannian steepest descent method. 

\section{Efficient approximation of geodesic equations}
\label{sec:geodesicApproximation}
This section aims to survey possible techniques that can be used to numerically solve the geodesic equation given in~\cite{Michor2007} or efficiently approximate it. 

\paragraph{Direct approach} as suggested in \cite{HerzogLoayzaRomero:2020:1}, one can use the Hamiltonian formulation of the geodesic equation and use a symplectic numerical integrator to solve the associated Hamiltonian system. A prominent example of such methods is the St\"ormer-Verlet scheme; see, e.g.,~\cite{hairer2003geometric}. 

\paragraph{Symplectic model order reduction} when dealing with large systems of PDEs there is always the option to use model order reduction (MOR). However, to preserve the structure of the Hamiltonian system one need to focus on symplectic model order reduction as the one developed in~\cite{peng2016symplectic}. 

\paragraph{Hamiltonian sparse identification of nonlinear dynamics} another option is to use the deep learning technique of sparse identification of nonlinear dynamics (SINDy) (cf.~\cite{brunton2022data}) and apply it to Hamiltonian system. This has not been used in the context of geodesic equations but it has been applied to standard Hamiltonian systems; see, e.g., \cite{canizares2024symplectic}.

\paragraph{Discrete geodesic calculus} additionally, it is also possible to exploit the nonlinear structure of the diffeomorphism group and develop the so-called discrete geodesic calculus, introduced in \cite{rumpf2013discrete}. The main idea of this approach is to locally approximate the Riemannian distance by a computationally inexpensive dissimilarity measure. Based on this dissimilarity measure, one can develop a strategy to approximate discrete geodesic, and thus the exponential map from a recursive expression. 

\section{Conclusion}
The use of geodesics instead of retractions in optimization procedures creates opportunities for implementing higher-order optimization methods. To formulate these methods on Riemannian shape spaces, we need to incorporate specific operators: retractions for first-order methods and vector transports for higher-order approaches. These operators are determined by the chosen Riemannian metric and often involve complex calculations.
The discrete geodesic calculus approach described by Rumpf and Wirth in~\cite{rumpf2013discrete} is particularly significant in this context. By defining a dissimilarity measure, this method approximates not only exponential maps but also logarithmic functions and---crucially---parallel transport. This parallel transport capability is essential for developing Quasi-Newton methods on Riemannian manifolds. Furthermore, this approach defines a discrete Riemannian connection that can be used to compute the Riemannian Hessian, enabling Newton iterations for shape optimization problem on shape manifolds.

 \bibliographystyle{splncs04}
 \bibliography{literature}

\end{document}